\documentclass[12pt]{article}
\usepackage{amsmath,amssymb,amsbsy,amsfonts,amsthm,
latexsym,amsopn,amstext,amsxtra,euscript,amscd}
\usepackage{amsmath,amsthm,epsfig}

\newcommand{\Q}{{\mathbb Q}}
\newcommand{\Z}{{\mathbb Z}}

\newcommand{\D}{{\mathbb D}}

\newcommand{\LL}{{\mathbb L}}
\newcommand{\CC}{{\mathcal C}}

\newcommand{\beq}{\begin{equation}}
\title{\bf Heron triangles with two fixed sides }
\author{
{\sc Eugen J. Ionascu}\\
{\rm Columbus State University}\\
{4225 University Avenue}\\
{Columbus, GA 31907}\\
{\tt ionascu\_eugen@colstate.edu} \vspace{0.5cm}\\
{\sc Florian~Luca}\\
{\rm Instituto de Matem{\'a}ticas}\\
{ Universidad Nacional Aut\'onoma de M{\'e}xico} \\
{C.P. 58180, Morelia, Michoac{\'a}n, M{\'e}xico} \\
{\tt fluca@matmor.unam.mx} \vspace{0.5cm}\\
{\sc Pantelimon St\u anic\u a}\\
{Auburn University Montgomery}\\
{Department of Mathematics}\\
{Montgomery, AL 36124-4023, USA}\\
{\tt pstanica@mail.aum.edu}
}
\begin{document}

\date{\today}

\maketitle


\begin{abstract}
In this paper, we study the function $H(a,b)$, which associates to every pair
of positive integers $a$ and $b$ the number of positive integers $c$ such that
the triangle of sides $a,b$ and $c$ is Heron, i.e., has integral area.
In particular, we prove that $H(p,q)\le 5$ if $p$ and $q$ are primes, and that
$H(a,b)=0$ for a random choice of positive integers $a$ and $b$.
\end{abstract}


\def\RR{{\rm I}\!{\rm R}}
\def\fp#1{(#1)}

\newtheorem{thm}{Theorem}
\newtheorem{rem}{Remark}
\newtheorem{lem}{Lemma}[section]
\newtheorem{theorem}[lem]{Theorem}
\newtheorem{alg}[lem]{Algorithm}
\newtheorem{cor}[lem]{Corollary}
\newtheorem{conj}[lem]{Conjecture}
\newtheorem{prop}[lem]{Proposition}
\newtheorem{heu}[lem]{Heuristic}
\newtheorem{LMN}[lem]{Theorem}
\newtheorem{exa}{Example}

\def\D{{\cal D}}
\def\RR{{\rm I}\!{\rm R}}
\def\LL{{\rm I}\!{\rm L}}
\def\NN{{\rm I}\!{\rm N}}
\def\MM{{\rm I}\!{\rm M}}
\def\QQ{{\rm I}\!\!\!{\rm Q}}
\def\ZZ{{\rm Z}\!\!{\rm Z}}
\def\CC{{\rm I}\!\!\!{\rm C}}
\def\vp{\varepsilon}
\def\phi{\varphi}
\def\ra{\rightarrow}
\def\sd{\bigtriangledown}
\def\ac{\mathaccent94}
\def\wi{\sim}
\def\wt{\widetilde}
\def\bb#1{{\Bbb#1}}
\def\bs{\backslash}
\def\cal{\mathcal}
\def\ca#1{{\cal#1}}
\def\Bbb#1{\bf#1}
\def\blacksquare{{\ \vrule height7pt width7pt depth0pt}}
\def\bsq{\blacksquare}
\def\n{\noindent}
\def\iter#1{^{(#1)}}
\def\bh{\ca B(H)}
\def\ld{\overline}
\def\vpe{V_{\psi}^{\eta}}
\def\wu{\ca W(\ca U)}
\def \dper {\ca D\ca P}
\def\lir{L^{\infty}(\RR)} \def\ltr{L^2(\RR)}
\def\ltrn{L^2(\RR^n)}
\def\lie{L^{\infty}(E)}\def\lte{L^2(E)}
\def\mpe{\ca M_{\psi,\eta}}
\def\wh{\widehat}
\def\eproof{$\hfill\bsq$\par}
\def\ws{\ca W\ca S}
\def\wi{\ca W\ca I}
\def\dst{\displaystyle}
\def\wsna{{\ca W}{\ca S}(n,A) }
\def\daw{{\ca W}_A}
\def\fwi{\ca W\ca I_1}
\def\swi{\ca W\ca I_2}
\def\ds{\displaystyle}
\def\du{\overset{\text {\bf .}}{\cup}}
\def\Du{\overset{\text {\bf .}}{\bigcup}}
\def\b{$\blacklozenge$}

\let\workingver=n
\def\begeq#1{\begin{equation}\mylabel{#1}}
\def\endeq{\end{equation}}
\def\begalg{\begin{alg}}
\def\endalg{\end{alg}}
\def\be{\beta}
\def\be{\beta}
\def\de{\delta}
\def\st{\star}
\def\op{\oplus}
\def\Llr{\Longleftrightarrow}
\def\Om{\Omega_f}
\def\rs{{\it RotS} }
\def\hF{{\hat {\cal F}}}
\def\refeq#1{\if\workingver y(\ref{#1})-
[[#1]]\else(\ref{#1})\fi}
\def\refth#1{\if\workingver y\ref{#1}-[[#1]]\else\ref{#1}\fi}
\def\mylabel#1{\if\workingver y\label{#1}{\bf\ \ [[#1]]\ \ }
\else\label{#1}\fi}
\def\mybibitem#1{\if\workingver y\bibitem{#1}{\bf\ \ [[#1]]\ \
}
\else\bibitem{#1}\fi}
\newcommand{\kro}[2]{\left( \frac{#1}{#2} \right) }

\section{Introduction}
\label{intro}

A {\em Heron triangle} is a triangle having the property that the
lengths of its sides as well as its area are positive integers.
There are several open questions concerning the existence
of Heron triangles with certain properties. For example
(see \cite{Guy94}, Problem D21), it is not known whether there
exist Heron
triangles having the property that the lengths of all of their
medians are positive integers, and it is not known (see
\cite{Harborth-Kemnitz} and \cite{Harborth et al. 96})
whether there exist Heron triangles having the property
that the lengths of all of their sides are Fibonacci numbers
other than
the triangle of sides $(5,5,8)$.
A different unsolved problem which asks for the existence
of a perfect cuboid; i.e., a rectangular box having the
lengths of all the sides, face diagonals, and main diagonal
integers has been related (see \cite{Luca 00}) to the
existence
of a Heron triangle having the lengths of its sides perfect
squares and the lengths of its bisectors positive integers
(see \cite{CG05} for more recent work).


Throughout this paper, we use the Vinogradov symbols $\gg$ and $\ll$ and
the Landau symbols $O$ and $o$ with their usual meanings. 
We write $\log x$ for the maximum between the natural logarithm of $x$
and $1$. Thus, all logarithms which will appear are $\ge 1$.

\section{Heron triangles with two fixed prime sides}

In this paper, we investigate (motivated by \cite{MAA04}) the
Heron triangles
with two fixed sides. Given a Heron triangle of sides
$a,~b$ and $c$ we write $S$ for its area and $s=(a+b+c)/2$
for its semiperimeter. It is known that $s$ is an integer and
that
$\min\{a,b,c\}\ge 3$. It is also known that if the triangle is
isosceles,
say if $a=b$, then $c$ is even, and if we write $h_c$ for
the length
of the altitude which is perpendicular on the side $c$, then
$h_c$ is an integer (see, for example, \cite{Luca 03}). In
particular,
$(c/2,h_c,a)$ is a Pythagorean triple with $a$ as hypothenuse.

\medskip

Now let $a$ and $b$ be positive integers
and let $H(a,b)$ be the number of distinct values of $c$, such
that $a,b$ and $c$ are the sides of a Heron triangle. In this
paper,
we investigate the function $H(a,b)$.

We start with an easy proposition that finds some bounds for $H(a,b)$.

\begin{prop}
\label{prop:1}
If $a\le b$ are fixed, then $0\le H(a,b)\le 2a-1.$
\end{prop}
\begin{proof}
Apply the triangle inequality to obtain
$ c+a>b$ and $ a+b>c,$
which imply that $b-a<c<b+a$. The inequality on $H(a,b)$
follows.
\end{proof}

Assume now that the prime power factorization of
$$
n=2^{a_0} p_1^{a_1}\cdots p_s^{a_s}q_1^{e_1}\cdots q_r^{e_r},
$$
where
$p_i\equiv 3\pmod 4$ for $i=1,\dots,s$ and
$q_j\equiv 1\pmod 4$ for $j=1,\dots,r$.
We let $d_1(n)=0$ if any of  the $a_j$ with $j\ge 1$ is odd.
Otherwise, let
$d_1(n)=\prod_{j=1}^r (1+e_j)$ be the number of
divisors of $n$ with primes among the $q_i$'s.
We shall need the following classical result
(see, for instance, \cite{Beiler66}).
\begin{lem}
\label{lem:1}
$(i) [Euler\ (1738)]$ A positive integer $n$ can be represented as
the sum of two
squares if and
only if $a_i\equiv 0\pmod 2$ for all $i=1,\ldots,s$.\\
$(ii)$ Assume that $a_i\equiv 0\pmod 2$ for all
$i=1,2,\ldots,s$.
The number of representations of $n$, say $r_2(n)$, as the sum
of two squares
ignoring order is given by $d_1(n)/2$ if $d_1(n)$
is even,
and $(d_1(n)-(-1)^{a_0})/2$ if $d_1(n)$ is odd.
\end{lem}

Using Lemma \ref{lem:1},
one can obtain a much better result than our previous
Proposition \ref{prop:1}.
Let $\tau(n)$ stand for the number of divisors of the positive
integer $n$.
\begin{theorem}
\label{thm:tau}
If $a$ and $b$ are fixed, then $H(a,b)\le 4\tau(ab)^2$.
\end{theorem}
\begin{proof} We start with the following observation.
Let $\gamma$ be any angle of a Heron triangle. Suppose that
$\gamma$ opposes the side $c$. Both $\cos \gamma$ and
$\sin \gamma$ are rational,
since $\cos\gamma=(a^2+b^2-c^2)/(2ab)\in \Q$ and
$\sin\gamma=2S/(ab)\in
\Q$.
Thus, if $\sin\gamma=u/v$, then $v=m^2+n^2$
for some coprime integers $m$ and
$n$ of opposite parities and $u\in \{2mn,~|m^2-n^2|\}$.
Since $a$ and $b$ are fixed,
in order to determine the number of Heron triangles that can
be constructed
using $a$ and $b$, it suffices to determine the number of
possible
angles $\gamma$ between $a$ and $b$. Since $ab\sin
\gamma=2S\in \Z$, if $\sin \gamma=u/v$, then
$v|ab$. For every divisor $v$ of $ab$, the number of possible
values for
$\sin\gamma$ is given by twice the number of ways of writing
$v =m^2+n^2$ with coprime positive integers
$m$ and $n$, because for every such representation we may have
$u=|m^2-n^2|$ or $u=2mn$. Furthermore, since when $\sin \gamma$ is fixed
there are
only two possibilities for $\gamma$, it follows that if $v|ab$
is fixed,
then the number of values for $\gamma$ is, by Lemma
\ref{lem:1},
at most $4d_1(v)$. Hence,
\[
H(a,b)\le 4\sum_{v|ab} d_1(v)\le 4\tau(ab)\sum_{v|ab} 1\le
4\tau(ab)^2.
\]
\end{proof}
When $a$ and $b$ are prime numbers, we
obtain a more precise result which improves upon \cite{MAA04}.

\begin{theorem}
\label{thm-main1}
If  $p$ and $q$ are two fixed primes, then
$$
H(p,q)\ \text{is} \
\begin{cases}
=0\ \  \text{if both p and}\ q\ \text{are}\ \equiv 3 \pmod 4,\\
=2\ \  \text{if}\ p=q\equiv 1 \pmod 4,\\
\le 2\ \text{if}\  p\not=q\  \text{and exactly one of p and q is
}\  \equiv 3
\pmod 4, \\
\le 5\ \text{if}\  p\not=q\  \text{and both p and q are } \equiv
1 \pmod 4.
\end{cases}
$$
\end{theorem}
\begin{proof}   Consider a Heron triangle of sides $p,~q$ and $x$.
Let $s=(p+q-x)/2$ be its semiperimeter.
By the Heron formula, its area is
\begin{equation}
\label{heron-formula}
S = \sqrt{s(s - p)(s - q)(s - x)},
\end{equation}
which after squaring becomes
\begin{equation}
\label{heron-equiv}
16S^2 =2p^2q^2 + 2p^2x^2 + 2q^2x^2 -p^4 - q^4 - x^4.
\end{equation}

Since $\min\{p,q,x\}\ge 3$ and $s$ is an integer, it follows
that both
$p$ and $q$ are odd primes and $x$ is an even integer.

Assume first that the triangle is isosceles. Since both $p$
and $q$ are odd
and $x$ is even, we get $p=q$. Furthermore, $(x/2,h_x,p)$
is a Pythagorean triple. This shows that $p=u^2+v^2$ with
coprime
positive integers $u$ and $v$ and $\{x/2,h_p\}=\{2uv,|u^2-
v^2|\}$. This leads
to either $x=4uv$, or $x=2|u^2-v^2|$. Hence, $H(p,p)\le 2$.
Furthermore,
$H(p,p)=0$ unless $p$ is a sum of two squares, which happens precisely
when $p\equiv 1\pmod 4$.

Assume now that $2<p<q$. Equation (\refth{heron-equiv}) is
equivalent to
\[
(p^2+q^2-x^2)^2+(4S)^2=(2pq)^2.
\]
With $x=2y$, the above equation can be simplified to
\begin{equation}\label{reduced}
((p^2+q^2)/2-2y^2)^2+(2S)^2=p^2q^2.
\end{equation}

Let us notice that if $(p^2+q^2)/2-2y^2=0$, we then get
a right triangle with legs $p$ and $q$ and hypothenuse $2y$,
which
does not exist because both $p$ and $q$ are odd.
Thus, we are lead to considering only the representations of
$p^2q^2$
as a sum of two positive squares.
We know, by Lemma \ref{lem:1} again, that the number of such
representations
(disregarding order) is
$0$ if both $p$ and $q$ are $\equiv 3\pmod 4$, it
is $(d_1(p^2q^2)-1)/2=1$
if exactly one of $p$ and $q$ is $\equiv 3\pmod 4$,
and finally it is
$(d_1(p^2q^2)-1)/2=4$ if neither $p$ nor $ q$ is $\equiv 3\pmod 4$.
Since in the representation (\ref{reduced}), the positive integer $2S$ is even
(thus, $(p^2+q^2)/2-2y^2$ is odd), each such representation
generates
only two possible alternatives.
When $p$ (or $q$) is congruent to $1$ modulo $4$,
we then write $p=u^2+v^2$ (or $q=z^2+w^2$) with $u$
(respectively $z$) odd.
The representations of $p^2q^2$ as $\left((p^2+q^2)/2-
2y^2\right)^2+(2S)^2$
are
\begin{eqnarray}
&&\left(p(z^2-w^2)\right)^2+(2pzw)^2,\ \label{eq4}\\
&&\left(q(u^2-v^2)\right)^2+(2quv)^2,\ \label{eq5}\\
&&\left((u^2-v^2)(z^2-w^2)+ 4uvzw\right)^2+
\left((u^2-v^2)2zw-(z^2-w^2)2uv\right)^2,\ \label{eq6}\\
&& \left((u^2-v^2)(z^2-w^2)- 4uvzw\right)^2+
\left((u^2-v^2)2zw+ (z^2-w^2)2uv\right)^2.\ \label{eq7}
\end{eqnarray}
In the case that both $p$ and
$q$ are $\equiv 3\pmod 4$,
we have no nontrivial representation of $p^2q^2$
as a sum of two squares, therefore $H(p,q)=0$.
If only one of these two primes is congruent to $1$ modulo $4$,
then there is only one representation, namely either
(\ref{eq4}) or
(\ref{eq5}), according to whether $q\equiv 1 \pmod 4$, or
$p\equiv 1\pmod 4$, which then leads to the pair equations
(\ref{eq8}) and
(\ref{eq9}), or (\ref{eq10}) and (\ref{eq11})
below, respectively.
Thus, in this case $H(p,q)\le 2$.
In the last and most interesting case when both
$p$ and $q$ are $\equiv 1\pmod 4$, each of the
equations below
could (at least at a first glance) produce an integer solution $x_i$ (here,
$x_i=2y$):
\begin{eqnarray}
&&p^2+q^2+2p(z^2-w^2)=x_1^2,\label{eq8}\\
&&p^2+q^2-2p(z^2-w^2)=x_2^2,\label{eq9}\\
&&p^2+q^2+2q(u^2-v^2)=x_3^2,\label{eq10}\\
&&p^2+q^2-2q(u^2-v^2)=x_4^2,\label{eq11}\\
&&p^2+q^2-2(u^2-v^2)(z^2-w^2)-8uvzw=x_5^2,\label{eq12}\\
&&p^2+q^2-2(u^2-v^2)(z^2-w^2)+8uvzw=x_6^2,\label{eq13}\\
&&p^2+q^2+2(u^2-v^2)(z^2-w^2)-8uvzw=x_7^2,\label{eq14}\\
&&p^2+q^2+2(u^2-v^2)(z^2-w^2)+8uvzw=x_8^2.\label{eq15}
\end{eqnarray}

To conclude the proof of Theorem \ref{thm-main1}, it suffices to show that
if $p$ and $q$ are fixed, then at most five of the above eight
equations can produce integer solutions $x_i$.

\medskip

We first study the case when
equation (\ref{eq8}) has an integer solution $x_1$.

\begin{lem}
\label{Lem2.3}
Assume that $q=w^2+z^2$ with positive integers $w$ and $z$, where $z$ is odd.
Then the prime $p$ and the positive integer $x_1$ satisfy (\ref{eq8}) if and
only if
$p=\frac{z^2}{\delta}+\frac{w^2}{\delta+1}$,
where $\delta$ is a positive integer such that $\delta$ divides $z^2$
and $\delta+1$ divides $w^2$.
\end{lem}

\begin{proof}  Let us prove the ``if" part. Substituting the above expressions
of $p$ and $q$ in terms of $z^2$ and $w^2$ into the equality (\ref{eq8}),
we get
\begin{eqnarray*}
&&p^2+q^2+2p(z^2-
w^2)=z^4\left(\frac{1}{\delta^2}+\frac{2}{\delta}+1\right)+
w^4\left(\frac{1}{(\delta+1)^2}-
\frac{2}{(\delta+1)}+1\right)\\
&&+z^2w^2\left(\frac{2}{\delta(\delta+1)}+2+\frac{2}{\delta+1}
-
\frac{2}{\delta}\right)=
z^4\left(\frac{1+\delta}{\delta}\right)^2+
w^4\left(\frac{\delta}{\delta+1}\right)^2+2z^2w^2 \\
&&=\left(\frac{(\delta+1)z^2}{\delta}+\frac{\delta
w^2}{\delta+1}
\right)^2.
\end{eqnarray*}
So, $\displaystyle x_1=\frac{(\delta+1)z^2}{\delta}+
\frac{\delta w^2}{\delta+1}$,
which is an integer under our assumptions.\par

For the other implication,
rewrite (\ref{eq8}) as $(p+z^2-w^2)^2+(2zw)^2=x_1^2$.
Using the general form of Pythagorean triples,
we may write
\begin{eqnarray}
&& p+z^2-w^2=d(m^2-n^2)\ \text{and}\
2zw=2dmn\label{neweq16},~\text{\rm or}\\
&&p+z^2-w^2=2dmn\ \text{and}\ 2zw=d(m^2-n^2),\label{neweq17}
\end{eqnarray}
for some nonzero integers $m,n,d>0$, with $m$ and $n$
coprime and of
opposite parities.

\par
We first look at the instance \refeq{neweq16}.
Since $dmn=zw>0$, we may assume that both $m$ and $n$ are
positive integers. Solving the above equation for $m$ and
substituting it into $p+z^2-w^2=d(m^2-n^2)$ gives
$$
p=w^2-z^2+dm^2-dn^2=\frac{(z^2+dn^2)(w^2-dn^2)}{dn^2}.
$$
Since $p$ is prime and $z^2+dn^2>dn^2$, we see that it is
necessary that
$dn^2=(w^2-dn^2)\delta$ and $z^2+dn^2=\delta p$ for some
positive integer $\delta$.
The first relation can be equivalently written as
$w^2\delta=(\delta+1)dn^2$. Because
$\delta$ and $\delta+1$ are coprime, it follows that $\delta+1$ divides
$w^2$.

\medskip

Similarly, solving for $n$ in \refeq{neweq16} and arguing as
before,
we obtain $\delta'(dm^2-z^2)=dm^2$ and $w^2+dm^2=\delta' p$
for some positive integer $\delta'$. {From} these relations,
we get $(\delta+1)(\delta'-1)=
\delta\delta'$, or $\delta'=\delta+1$.
Hence, $z^2(\delta+1)=\delta dm^2$, which implies that $\delta$
divides $z^2$.
We then get
$$p=w^2-z^2+dm^2-dn^2=w^2-z^2+\frac{z^2(\delta+1)}{\delta}-
\frac{w^2\delta}{(\delta+1)}=\frac{z^2}{\delta}+\frac{w^2}{\delta+1},
$$
which completes the proof of the lemma in this case.\par

We now look at the instance \refeq{neweq17}.
If we denote by $s=m+n$ and $t=m-n$, then $p+z^2-w^2=d(s^2-
t^2)/2$,
and $2zw=dst$. Note that $s$ and $t$ are both odd and, in
particular,
$st\not=0$ and $d$ is even.
As before, solving for $s$, we get
$s=\displaystyle\frac{2zw}{dt}$,
and substituting it into $p+z^2-w^2=d(s^2-t^2)/2$ we get
$\displaystyle p=\frac{(2z^2+dt^2)(2w^2-dt^2)}{2dt^2}$,
which implies, as before and since
$2z^2+dt^2>dt^2$, that both $\mu (2w^2-dt^2)=2dt^2 $ and
$2z^2+dt^2=p\mu $
hold with some positive integer $\mu$. Hence,
$2w^2\mu=(\mu+2)dt^2$.

\medskip

\par
As before, solving now for $t$ we get $\displaystyle
t=\frac{2zw}{ds}$,
and substituting it into $p+z^2-w^2=d(s^2-t^2)/2$,
we arrive at $\displaystyle p=\frac{(ds^2-
2z^2)(2w^2+ds^2)}{2ds^2}$.
By a similar argument, we get $2w^2+ds^2=p\mu'$
and $2ds^2=\mu'(ds^2-2z^2)$ with some positive integer $\mu'$.
This implies that $(\mu'-
2)ds^2=2\mu'z^2$.
{From} these relations, we obtain $(\mu'-2)(\mu+2)(dst)^2=\mu
\mu' 4z^2w^2$.
This shows that $\mu'=\mu+2$. Since $z$ and $w$ are of
opposite parities,
we get that $\mu$ is even. The conclusion of the lemma
is
now verified by taking $\delta=\mu/2$.
\end{proof}
\begin{rem}
Similar statements hold for relations
$\refeq{eq9}$--$\refeq{eq11}$, as well.
\end{rem}

We now give a numerical example.

\begin{exa}
Let $s$ and $t$ be positive integers such that both $p=3t^2+2s^2$ and
$q=9t^2+4s^2$ are prime numbers. If we write $x_1=6(s^2+t^2)$,
then the triangle with sides $(p,q,x_1)$ is Heron,
and $x_1$ satisfies equation (\ref{eq8}).
\end{exa}
Another interesting example is $p=5$ and $q=1213=27^2+22^2$,
which gives $z=27$ and $w=22$.
In this case equation (\ref{eq9})
is satisfied with $x_2=1212$, which leads to the Heron triangle
with sides $(5,1213,1212)$.
   We note that $\frac{27^2}{243}+\frac{22^2}{242}=3+2=p$,
which shows that $\delta=242$.
   This example shows that $\delta$  does
not have to divide $z$ or $w$. \par
We now record the following corollary to Lemma \ref{Lem2.3}.

\begin{cor}
\label{cor:1}
Assume that $q=z^2+w^2$ is an odd prime. Then
equations (\ref{eq8}) and (\ref{eq9})
cannot both have integer solutions $x_1$ and
$x_2$.
\end{cor}
\begin{proof}  Assume, by way of contradiction, that both
$x_1$ and $x_2$ are integers.
Then, by Lemma \refth{Lem2.3}, we must have
$$
p=\frac{z^2}{\delta}+\frac{w^2}{\delta+1}=\frac{z^2}{\mu+1}+
\frac{w^2}{\mu}
$$
for some positive integers $\delta$ and $\mu$
satisfying the conditions from the statement of Lemma
\refth{Lem2.3}.
If $\mu\ge \delta+1$, then
$$
\frac{z^2}{\mu+1}+\frac{w^2}{\mu}<\frac{z^2}{\delta}+\frac{w^2
}{\delta+1}.
$$
Similarly,  if $\delta \ge \mu+1 $, we then have
$$
\frac{z^2}{\mu+1}+\frac{w^2}{\mu}>\frac{z^2}{\delta}+\frac{w^2
}{\delta+1}.
$$
Thus, we must have $\delta=\mu$.
But this is impossible because $z$ and $w$ have opposite
parity.
\end{proof}

\begin{cor}
Assuming that both $p=u^2+v^2$ and $q=z^2+w^2$ are
odd primes, then among the four equations (\ref{eq8}),
(\ref{eq9}),
   (\ref{eq10}) and (\ref{eq11}), at most one of the $x_i$'s is
an integer
($i=1,\dots,4$).
\end{cor}
\begin{proof}  By way of contradiction, suppose that more than one
of the
$x_i$'s is an integer. By Corollary \ref{cor:1}, the number
of
such cannot exceed two. We now use Lemma \refth{Lem2.3}
and obtain an integer $\delta$ with the property that either
$$
p= \frac{z^2}{\delta}+\frac{w^2}{\delta+1}\qquad \text{\rm
or}\qquad
p=\frac{w^2}{\delta}+\frac{z^2}{\delta+1}.
$$
In either situation, we get $p<q$.
Similarly, because one of equations (\ref{eq10}) or
(\ref{eq11})
must have an integer solution, we must have either
$$
q=\frac{u^2}{\delta'}+\frac{v^2}{\delta'+1}\qquad \text{\rm
or}\qquad
q=\frac{v^2}{\delta'}+\frac{u^2}{\delta'+1},
$$
which leads to $q<p$. This contradiction shows that at most
one of
the numbers $x_i$, $i=1\ldots 4$, is an integer.
\end{proof}

This completes the proof of our Theorem \refth{thm-main1}.
\end{proof}

\medskip

We now give an example of an instance in which both
equations
(\ref{eq12}) and (\ref{eq14}) have integer solutions.

\begin{prop}
Both equations (\ref{eq12}) and (\ref{eq14}) have integer
solutions
if $p=(ij)^2+(kl)^2$ and $q=(ik)^2+(jl)^2$
for some positive integers $i,j,k$ and $l$ with $i,~j,~k$ odd and
$l$ even.
\end{prop}
\begin{proof} Since $u$ and $z$ are odd, we have $u=ij$, $v=kl$,
$z=ik$ and $w=jl$. A straightforward calculation shows that
\[
\begin{array}{l}
p^2+q^2-2(u^2-v^2)(z^2-w^2)-
8uvzw=(i^2j^2+k^2l^2)^2+(i^2k^2+j^2l^2)^2-\\
-2(i^2j^2-k^2l^2)(i^2k^2-j^2l^2)-8i^2j^2k^2l^2=(i^2j^2-
k^2l^2)^2+\\
+(i^2k^2-j^2l^2)^2-2(i^2j^2-k^2l^2)(i^2k^2-j^2l^2)=
\left[(i^2j^2-k^2l^2-(i^2k^2-j^2l^2)\right]^2=\\
(k^2-j^2)^2(i^2+l^2)^2=x_5^2,
\end{array}
\]
with $x_5=|k^2-j^2|(i^2+l^2)$.

Similarly,
\[
\begin{array}{l}
p^2+q^2+2(u^2-v^2)(z^2-w^2)-
8uvzw=(i^2j^2+k^2l^2)^2+(i^2k^2+j^2l^2)^2+\\
+2(i^2j^2-k^2l^2)(i^2k^2-j^2l^2)-8i^2j^2k^2l^2=(i^2j^2-
k^2l^2)^2+\\
+(i^2k^2-j^2l^2)^2+2(i^2j^2-k^2l^2)(i^2k^2-j^2l^2)=
\left[(i^2j^2-k^2l^2+(i^2k^2-j^2l^2)\right]^2=\\
(k^2-j^2)^2(i^2+l^2)^2=x_7^2,
\end{array}
\]
with $x_7=|i^2-l^2|(j^2+k^2)$.
\end{proof}

The examples below are of the type, which additionally
satisfy three of the equations (\ref{eq12})--(\ref{eq15}):

$p=21521, q= 14969, {x_1=14952, x_2=15990, x_3=33448}$,

$p=4241, q=2729, {x_1=1530, x_2=1850, x_3=6888}$,

$p=898361, q=161009, {x_1=952648, x_2=896952, x_3=870870}$,

$p=659137, q=252913, {x_1=512720, x_3=688976, x_3=722610}$,

$p=4577449, q=11893681, {x_1=11843832, x_2= 14876232,
x_3=10174630}.$
We currently have no example satisfying all equations
(\ref{eq12})--(\ref{eq15}), in the same time.

 Let us observe that the equations (\ref{eq12})--(\ref{eq15})
are equivalent in the written order to
\begin{equation}
(u^2-v^2-z^2+w^2)^2+4(uv-zw)^2=x_5^2,
\end{equation}
\begin{equation}
(u^2-v^2-z^2+w^2)^2+4(uv+zw)^2=x_6^2,
\end{equation}
\begin{equation}
(u^2-v^2+z^2-w^2)^2+4(uv-zw)^2=x_7^2,
\end{equation}
\begin{equation}\label{last}
(u^2-v^2+z^2-w^2)^2+4(uv+zw)^2=x_8^2.
\end{equation}
In particular, we have $H(p,q)\ge 2$ if $uv=zw$, or $u^2-
v^2=z^2-w^2$ or $u^2-v^2=w^2-z^2$.

\begin{theorem}\label{second}
Assume that $q=z^2+w^2$ is an odd prime, $u$, $v$, $z$ and $w$
are positive integers with both $u$ and $z$ odd and both $v$ and $w$
even. Then equation (\ref{last}) has integer solutions
if and only if
there exist positive integers $m,n,a$ and $b$ with $m$ and $n$ coprime
of opposite parities and $a$ and $b$ also coprime,
such that $a$ divides $nw-mz$, $b$ divides $nz+mw$,
and
\begin{equation}\label{uandv}
\begin{cases} u=nk-\frac{zb}{a} \in \Z \\ v=mk-
\frac{wb}{a} \in \Z,
\end{cases}
\end{equation}
where $k=\frac{(a^2+b^2)(nz+mw)}{ab(n^2+m^2)}$.
\end{theorem}
\begin{proof} For the sufficiency part, let us first calculate $u^2-
v^2+z^2-w^2$. We have, from (\ref{uandv}),
\[
\begin{array}{l}
\ds u^2-v^2+z^2-w^2=(n^2-m^2)k^2-2\frac{kb}{a}(nz-
mw)+\frac{b^2}{a^2}(z^2-w^2)+z^2-w^2=
\\ \\
\ds (n^2-m^2)k^2-2\frac{(a^2+b^2)(n^2z^2-
m^2w^2)}{a^2(n^2+m^2)}+\frac{(a^2+b^2)(z^2-w^2)}{a^2}.
\end{array}
\]
This reduces to
\[
\begin{array}{l}
\ds u^2-v^2+z^2-w^2= \ds (n^2-m^2)k^2-\frac{(a^2+b^2)(n^2-
m^2)(z^2+w^2)}{a^2(n^2+m^2)}=
\\ (n^2-m^2)\left(k^2-
\frac{(a^2+b^2)(z^2+w^2)}{a^2(n^2+m^2)}\right).
\end{array}
\]
On the other hand, $uv+zw=mnk^2-
\frac{kb}{a}(nw+mz)+\frac{b^2}{a^2}zw+zw$, which, as before,
gives
\[
\begin{array}{l}
\ds uv+zw=mnk^2-\frac{(a^2+b^2)[(nz+mw)(nw+mz)-
zw(n^2+m^2)]}{a^2(m^2+n^2)}=\\ \\
\ds =mn\left(k^2-
\frac{(a^2+b^2)(z^2+w^2)}{a^2(n^2+m^2)}\right).
\end{array}
\]
Putting these calculations together, we see that

$$(u^2-v^2+z^2-w^2)^2+4(uv+zw)^2=(m^2+n^2)^2\left(k^2-
\frac{(a^2+b^2)(z^2+w^2)}{a^2(n^2+m^2)}\right)^2.$$

\n This means that we can  take $x_8$ in (\ref{last}) to be
$\ds (m^2+n^2)k^2- (a^2+b^2)q/a^2$, which
must be an integer since we assumed $u$ and $v$ were integers.

For the necessity part, let us observe that $(u^2-v^2+z^2-
w^2)/2\equiv 1 \pmod 2$ and $uv+zw \equiv 0 \pmod  2$,
which imply that $u^2-v^2+z^2-w^2=d(n^2-m^2)$ and
$uv+zw=2dmn$ for some nonzero integers $m,n$ and $d>0$, with $m$
and $n$ coprime and of opposite parities. From here, we see
that $mn(u^2-v^2+z^2-w^2)-(n^2-m^2)(uv+zw)=0$.
This is equivalent to
\[
mn(u^2-v^2)-(n^2-m^2)uv+mn(z^2-w^2)-(n^2-m^2)zw=0,
\]
   or
\begin{equation}
\label{eqeq}
(nu+mv)(mu-nv)+(nz+mw)(mz-nw)=0.
\end{equation}
We are now ready to define our coprime numbers $a$ and $b$ by the
relation
\[
\frac{a}{b}=\frac{nu+mv}{nz+mw}.
\]
>From here, we see that $nz+mw=bs$ for some integer $s$, and
in light of (\ref{eqeq})
we also have
\[
\frac{a}{b}=\frac{nw-mz}{mu-nv},
\]
which shows that $nw-mz=as'$ for some integer $s'$. We are now
going to show that both equalities in (\ref{uandv}) hold.
Indeed, let us look at the last two equalites as a system of linear
equations with two unknowns in $u$ and $v$:
\begin{equation}
\begin{cases} nu+mv=\frac{a}{b} (nz+mw),\\
mu-nv=\frac{b}{a}(nw-mz).
\end{cases}\end{equation}
If one solves the above system for $u$ and $v$, one gets exactly the two
expressions in (\ref{uandv}). \end{proof}

{\bf Example:} If  $p=113=7^2+8^2$, $q=257=1^2+16^2$, which
means that $u=7$, $v=8$, $z=1$ and $w=16$. Because
$\frac{u^2-v^2+z^2-w^2}{uv+zw}=-15/4$, we get $n=1$ and $m=4$.
Then $\frac{nu+mv}{nz+mw}=\frac{3}{5}$, which gives $a=3$ and $b=5$.
Thus, $k=\frac{26}{3}$ and the formulae (\ref{uandv}) are
easily verified.

\begin{theorem}\label{important}
If $p=u^2+v^2$ and $q=z^2+w^2$ are odd distinct primes,  with
$u,v,z$ and $w$ as in Theorem ~ \ref{second} for which (\ref{last})
has integer solution $x_8$, then
\begin{equation}\label{pplusq}
p+q=\ell(s^2+{s'}^2),
\end{equation}
where $nz+mw=sb$, $nw-mz=s'a$ and $a^2+b^2=\ell (m^2+n^2)$ with
integers $s,~s'$ and $\ell$.
\end{theorem}

\begin{proof} We apply Theorem ~ \ref{second} and obtain the
numbers $a$, $b$, $m$, $n$ and $k$ with the specified
properties. Using the above notations, we get $\ds k=\frac{\ell
s}{a}$.  Squaring and adding together the equalities for $u$
and $v$ from Theorem ~ \ref{second}, we get
\[
p=u^2+v^2=(m^2+n^2)k^2-
2\frac{kb}{a}(nz+mw)+(z^2+w^2)\frac{b^2}{a^2},
\]
which becomes
\[
p=\frac{(m^2+n^2)\ell ^2s^2}{a^2}-2\frac{b^2\ell
s^2}{a^2}+q\frac{b^2}{a^2}=
\frac{(a^2+b^2)\ell s^2}{a^2}-2\frac{b^2\ell
s^2}{a^2}+q\frac{b^2}{a^2}.
\]
Equivalently,
\[
\begin{array}l
\ds p+q=\ell s^2+q\frac{a^2+b^2}{a^2}-\frac{b^2\ell
s^2}{a^2}=\ell s^2+\ell \frac{(z^2+w^2)(m^2+n^2)-
(nz+mw)^2}{a^2}=\\ \\
\ds =\ell s^2+\ell \frac{(nw-mz)^2}{a^2}=\ell(s^2+{s'}^2).
\end{array}
\]

Let us show now that $\ell$ is an integer. Since we know that
$u=nk-zb/a$ must be an integer, solving for $nak$ we obtain
that $n\frac{(a^2+b^2)s}{m^2+n^2}=nak=au+zb$ must be an integer as
well. Since $\gcd(m^2+n^2,n)=1$, we see that $m^2+n^2$
should divide $(a^2+b^2)s$. The conclusion we want follows
then provided that we can show that $m^2+n^2$ and $s$ are coprime.
To this end, it
is enough to prove that $\gcd(m^2+n^2,nz+mw)=1$.

By way of contradiction, assume that $m^2+n^2$ does not divide
$a^2+b^2$ and there exists a prime number $r$ which divides
both $m^2+n^2$ and $nz+mw$.
Then $r$ divides $n(nw-mz)=(m^2+n^2)w-m(nz+mw)$. Because
$\gcd(r,n)=1$, we see that $r$ must divide $nw-mz$.
Hence, $nq=n(z^2+w^2)=z(nz+mw)+w(nw-mz)$ is divisible by $r$.
The assumption that $q$ is prime together with the fact that
$\gcd(r,n)=1$ shows that $r=q$. In this case,
$m^2+n^2=(z^2+w^2)(m_1^2+n_1^2)$ with $n=m_1z-n_1w$ and
$m=n_1z+m_1w$
for some integers $m_1,n_1$. This gives $bs=nz+mw=m_1q$.

On the other hand, the fact that $m^2+n^2$ divides
$(a^2+b^2)(nz+mw)$ is equivalent to saying that $m_1^2+n_1^2$
divides $(a^2+b^2)m_1$. The fact that $\gcd(m_1^2+n_1^2,m_1)=1$
implies that
$m_1^2+n_1^2$ divides $a^2+b^2$. So, the fraction
$\frac{(a^2+b^2)(nz+mw)}{b(m^2+n^2)}$ can be simplified to
$\frac{a's}{q}$. If in addition $q$ divides $a'$, then
$m^2+n^2$ divides $a^2+b^2$, which we excluded. If $q$
divides $s$, then from the identity
$(m^2+n^2)(p+q)=(a^2+b^2)(s^2+{s'}^2)$ and the above observation,
we get $q(p+q)=a'(q^2s_1^2+{s'}^2)$. But this
implies that $q$ divides $s'$ as well, and after simplifying
both sides by $q$ we get that $q$ divides $p$,
which is a contradiction.\end{proof}

\begin{cor}\label{critical}
Under the assumptions of Theorem~\ref{important}, we have
\begin{equation}\label{alphabeta}
p+q=\alpha^2+\beta^2
\end{equation}
for some integers $\alpha$ and $\beta$.
\end{cor}
\begin{proof} We can apply Theorem ~ \ref{important}, as well as
the theory of Gaussian integers (or directly, by looking at the
prime power factorization), to get that
$\ell=\frac{a^2+b^2}{m^2+n^2}$ must be of the form
$m_1^2+n_1^2$ for some integers $m_1$ and $n_1$. The existence of
$\alpha$ and $\beta$ follows from the
Cauchy identity $(s^2+{s'}^2)(m_1^2+n_1^2)=(sm_1+s'n_1)^2+(sn_1-s'm_1)^2$.
\end{proof}

\section{Counting Heron triangles}

In this section, we prove the following theorem.
\begin{theorem}
\label{main-thm}
The estimate
\begin{equation}
\label{eq:main}
\sum_{a,b\le x} H(a,b)<x^{25/13+o(1)}
\end{equation}
holds as $x\to \infty$.
\end{theorem}

For the proof of the above result,
we will use the following completely explicit two dimensional
version of the Hilbert's Irreducibility Theorem due to
Schinzel and Zannier \cite{SZ}, whose proof is based on results
of Bombieri and Pila \cite{BP}.

\begin{lem}
\label{lem1}
Let $\Phi\in \Q[t,y]$ be a polynomial irreducible over $\Q$ of total
degree $D\ge 2$. Then, for every positive integer $\delta<D$ and for every
$T\ge 1$, the number of integer points $(t^*,y^*)$ such that $\Phi(t^*,y^*)=0$
and $\max\{|t^*|,|y^*|\}\le T$ is bounded by
$$
(3D\Delta)^{\Delta+4}T^{3/(3(\delta+3))},
$$
where $\Delta=(\delta+1)/(\delta+2)$.
\end{lem}
\begin{proof} This is Lemma 1 on page 294 in \cite{SZ}.
\end{proof}

\begin{lem}
\label{lem2}
Let $A,~B,~C,~D$ and $E$ be integers with $ADE\ne 0$. Then every irreducible
factor of
$$
F(t,z)=A t^4-(Bz^2+C)t^2-(Dz+E)^2\in \Q[t,z]
$$
has degree at least $2$ in $t$.
\end{lem}

\begin{proof} Assume that this is not so. Then there exists a factor of
$F(t,z)$ of the form $t-f(z)$ for some $f(z)\in \Q[z]$. Thus,
\begin{equation}
\label{eq:1}
Af(z)^4-(Bz^2+C)f(z)^2-(Dz+E)^2=0
\end{equation}
identically. By looking at the degrees,
it now follows easily that $f(z)$ is either linear or constant.
If $f(z)=K$ is constant, then the coefficient of $z$ in the above relation
is $-2DE\ne 0$, which is a contradiction.
Suppose now that $f(z)$ is linear. Since $f(z)\mid Dz+E$, we get that
$f(z)=\lambda(Dz+E)$ for some nonzero rational $\lambda$. Substituting
in \eqref{eq:1}, we get
$$
A\lambda^4(Dz+E)^2-\lambda^2(Bz^2+C)-1=0,
$$
which is impossible as the coefficient of $z$ is $2\lambda^4ADE\ne 0$.
\end{proof}

\begin{proof}[Proof of Theorem \refth{main-thm}]
We proceed in several steps.

\medskip

\noindent
  {\bf Notations.}
\par\noindent
We assume that $a,~b$ and $c$ are the sides of a Heron triangle $\Delta$
of surface
area $S$ opposing
the angles $\angle A,~\angle B$ and $\angle C$, respectively.
Since $ab\sin (\angle C)=2S\in \Z$ and $\cos C=(a^2+b^2-c^2)/(2ab)\in \Q$,
we get
that both $\sin C$ and $\cos C$ are rational. Furthermore, $\sin(C)=2S/(ab)$.
We let $a_{1,c}b_{1,c}$ be the denominator of the rational number
$\sin(\angle C)$ such that $a_{1,c}\mid a$ and $b_{1,c}\mid b$.
Since $a_{1,c}b_{1,c}$ is a sum of two coprime squares, it follows that
each one of $a_{1,c}$ and $b_{1,c}$ is also a sum of two coprime squares. We
then write
$$
a^2_{1,c}=U^2+V^2\qquad {\text{\rm and}}\qquad b^2_{1,c}=W^2+Z^2,
$$
with some integers $U,~V,~Z$ and $W$ such that
$\gcd(U,V)=\gcd(W,Z)=1$, $U\not\equiv
V\pmod 2$ and $Z\not\equiv W\pmod 2$. Furthermore, up to interchanging
$W$ and $Z$ and changing the signs of
some of $U,~V,~W$ and/or $Z$, we may assume that
$$\sin(\angle C)=\frac{UZ-VW}{a_{1,c}b_{1,c}}\qquad {\text{\rm and}}
\qquad
\cos(\angle C)=\frac{UW+VZ}{a_{1,c}b_{1,c}}.
$$Thus, we get that if we write $a_{2,c}=a/a_{1,c}$ and $b_{2,c}=b/b_{1,c}$,
then
\begin{eqnarray}
\label{eq:important}
c^2 &= & a^2+b^2-2ab\cos(\angle C)\nonumber\\
& = & a_{2,c}^2(U^2+V^2)+b_{2,c}^2(W^2+Z^2)-
2a_{2,c}b_{2,c}(UW+VZ)\nonumber\\
& = & (a_{2,c}U-b_{2,c}W)^2+(a_{2,c}V-b_{2,c}Z)^2.
\end{eqnarray}
Furthermore, we know that $a_{1,c}=u^2+v^2$, therefore
$\{U,V\}=\{|u^2-v^2|,2uv\}$ and similarly $b_{1,c}=w^2+z^2$, therefore
$\{W,Z\}=\{|w^2-z^2|,2wz\}$ for some integers $u,~v,~w$ and $z$ with
$\gcd(u,v)=1,~\gcd(w,z)=1$, $u\not\equiv v\pmod 2$ and
$w\not\equiv z\pmod 2$.

\medskip
\noindent
In a similar way, we get two more equations satisfied by our triangle
of the same shape as \eqref{eq:important}, obtained by circularly permuting
$a,~b$ and $c$ (and $\angle C$, $\angle A$ and $\angle B$, respectively).

\medskip
\noindent
{From} now on, we let $x$ be a large positive real number. We assume that
$\Delta$ is a Heron triangle of sides
$a,~b$ and $c$ satisfying $\max\{a,b,c\}\le x$.

\medskip
\noindent
{\bf Triangles having large $a_{2,b}$.}
\par\noindent Here, we assume that $\Delta$ is such that
\begin{equation}
\label{eq:min}
\min\{a_{2,b},~a_{2,c},~b_{2,a},~b_{2,c},~c_{2,a},~c_{2,b}\}\ge x^{4/5}.
\end{equation}
Then $a_{1,c}b_{1,c}=ab/(a_{2,c}b_{2,c})\le x^{2/5}$. Since
$a_{1,c}b_{1,c}=K^2+L^2$ for some coprime positive
integers $K$ and $L$ such that
$$
\sin(\angle C)\in \{2KL/(K^2+L^2),~|K^2-L^2|/(K^2+L^2)\},
$$
we get that $K<x^{1/5}$ and $L<x^{1/5}$. Thus, $\angle C$ can be chosen
in $O(x^{2/5})$ ways. Similarly, there are only $O(x^{2/5})$ choices for
$\angle B$. Since $a\le x$, and since a side and its two adjacent
angles determine the triangle uniquely, we get that the number of
possible triangles in this case is $O(x^{9/5})$.

\medskip
\noindent
{From} now on, we count only those Heron triangles $\Delta$ which do not
fulfill condition \eqref{eq:min}. To fix ideas, we assume that
$$
a_{2,c}=\min\{a_{2,b},~a_{2,c},~b_{2,a},~b_{2,c},~c_{2,a},~c_{2,b}\},$$
and that $a_{2,c}<x^{4/5}$. Furthermore,
to simplify notations, from now on we will drop
the subscript $c$. Relation \eqref{eq:important} becomes
\begin{equation}
\label{eq:pyt}
c^2=(a_{2}U-b_{2}W)^2+(a_{2}V-b_{2}Z)^2.
\end{equation}

\medskip
\noindent
  {\bf The case when $a_2V=b_2Z$.}
\par\noindent Fix $a_2$. Then $a_1=u^2+v^2\le x/a_2$, leading
to $|u|\le (x/a_2)^{1/2}$ and $|v|\le (x/a_2)^{1/2}$. Thus, there
are only $O(x/a_2)$ possibilities for the pair $(u,v)$. Since $V$ is
determined by the pair $(u,v)$, that is,
$V\in\{|u^2-v^2|,2uv \}$,
  it follows that there are only $O(x/a_2)$
possibilities for $V$. Now $a_2V=b_2Z$, therefore both $b_2$ and $Z$ divide
$a_2V$. If $Z=2wz$, then both $w$ and $z$ divide $a_2V$, while when
$Z=|w^2-z^2|$, then both $w-z$ and $w+z$ divide $a_2V$. In either case,
we get that the number of possibilities for the triple $(b_2,w,z)$
is $O(\tau(a_2V)^3)=x^{o(1)}$. It is now clear that $c$ is uniquely determined
once $a_2,~u,~v,~b_2,~w$ and $z$ are determined
as $c=|a_2U-b_2W|$. This argument shows that the number of possibilities
for $\Delta$ in this case is
$$
\le \sum_{a_2< x^{4/5}}\frac{x^{1+o(1)}}{a_2}=x^{1+o(1)}\sum_{a_2<x^{4/5}}
\frac{1}{a_2}=x^{1+o(1)}.
$$
A similar argument applies for the number of Heron triangles
$\Delta$ with $a_2U=b_2W$.

\medskip

\noindent
{From} now on, we assume that $a_2V\ne b_2Z$ and that $a_2U\ne b_2W$.
By relation \eqref{eq:pyt}, the three
numbers $c,~|a_2U-b_2W|$ and $|a_2V-b_2Z|$
form a Pythagorean triple with $c$ as hypothenuse. Thus,
there exist positive integers $d,~M$ and $N$ with $M$ and $N$ coprime and of
different parities such that
\begin{equation}
\label{eq:form}
\begin{split}
&c=d(M^2+N^2)\quad {\text{\rm and}}\\
&\{\,|a_2U-b_2W|,|a_2V-b_2Z|\,\}=
\{\,d|M^2-N^2|,2dMN\,\}.
\end{split}
\end{equation}

\medskip
\noindent
We shall only assume that
\begin{equation}
\label{eq:form1}
U=u^2-v^2,\qquad V=2uv,\qquad W=w^2-z^2,\qquad Z=2wz,
\end{equation}
and
\begin{equation}
\label{eq:form2}
a_2U-b_2W=d(M^2-N^2)\qquad {\text{\rm and}}\qquad a_2V-b_2Z=2dMN,
\end{equation}
since the remaing cases can be handled similarly.

\medskip
\noindent
We now take $\kappa_0$ to be a positive constant to be determined later.

\medskip
\noindent
  {\bf Triangles with large $d$.}
\par\noindent We assume that $d>x^{\kappa_0}$. Then $M^2+N^2<x^{1-\kappa_0}$,
therefore $|M|<x^{(1-\kappa_0)/2}$ and $|N|<x^{(1-\kappa_0)/2}$. Thus,
the pair $(M,N)$ can be chosen in $O(x^{1-\kappa_0})$ ways. We shall now fix
$a_2,~u$ and $v$. This triple can be fixed in $x^{1+o(1)}$ (first we fix
$a\le x$, then we fix $a_2\mid a$ and $u$ and $v$ such that $u^2+v^2\mid a$).
Dividing the last two equations \eqref{eq:form2} side by side and using
equation \eqref{eq:form1}, we get the equation
$$
\frac{a_2(u^2-v^2)-b_2(w^2-z^2)}{a_2uv-b_2wz}=\frac{M^2-N^2}{MN}.
$$
Cross multiplying and reducing modulo $b_2$, we get that
\begin{equation}
\label{eq:divi}
b_2\mid MNa_2(u^2-v^2)-a_2uv(M^2-N^2).
\end{equation}

\medskip
\noindent
We now distinguish the following cases.

\medskip
\noindent
{\bf Case 1. $MNa_2(u^2-v^2)-a_2uv(M^2-N^2)=0$.}
\par\noindent Then
$(u^2-v^2)/(uv)=(M^2-N^2)/(MN)$, and $u$ and $v$ are coprime. Thus, $uv=MN$,
giving only $x^{o(1)}$ possibilities for the pair $(u,v)$ once the pair
$(M,N)$ is fixed. Thus, if we first fix $c\le x$ and then
$d\mid c$ and $M^2+N^2\mid c$, we see that there are
$x^{1+o(1)}$ possibilities for the triple $(d,M,N)$. Further, once
$d,~M$ and $N$ are fixed, the pair $(u,v)$ can be chosen in
only $x^{o(1)}$ ways. Finally, $a_2$ is fixed
in $O(x^{4/5})$ ways. Thus, the number of ways to fix the pair $(a,c)$ in
this case is $x^{9/5+o(1)}$, and since $b$ can take only $x^{o(1)}$ ways
once the sides $a$ and $b$ are determined, we get that
the number of such Heron triangles $\Delta$ is $x^{9/5+o(1)}$.

\medskip
\noindent
  {\bf Case 2. $MNa_2(u^2-v^2)-a_2uv(M^2-N^2)\ne 0$.}
\par\noindent
Then $b_2$ can be chosen in $x^{o(1)}$ ways once $(a,M,N)$ are fixed.
Furthermore, we now get
\begin{equation}
\label{eq:quad}
b_2MN(w^2-z^2)-b_2(M^2-N^2)wz=MNa_2(u^2-v^2)-a_2uv(M^2-N^2).
\end{equation}
Thus, there are $O(x^{2-\kappa_0})$ choices for the triple
$(a,M,N)$ in this case, and then $b_2$ is fixed in only $x^{o(1)}$ ways.
Once $b_2,~d,~M$ and $N$ are fixed, equation \eqref{eq:quad} becomes
$$
g(w,z)=K,
$$
where $K$ is fixed and nonzero and $g(w,z)$ is the quadratic form
$$
Aw^2+Bwz+Cz^2,
$$
where
$$
A=b_2MN,\qquad B=-b_2(M^2-N^2)\qquad {\text{\rm and}}\qquad C=-b_2MN.
$$
The discriminant of $g$ is
$$
b_2^2((M^2-N^2)^2+4M^2N^2)=b_2^2(M^2+N^2)^2,
$$
which is a nonzero perfect square. It now follows easily that the resulting
equation \eqref{eq:quad} has only $x^{o(1)}$ possible solutions $(w,z)$.
Hence, the number of Heron triangles $\Delta$
is in this case $x^{2-\kappa_0+o(1)}$.

\medskip
\noindent
A similar count will show that the number of Heron triangles with
$b_2>x^{\kappa_0}$ is $x^{2-\kappa_0}$. Thus, we may assume that
both $b_2$ and $d$ are smaller than or equal to
$\kappa_0$. Since $a_2\le b_2$, we also
get that $a_2\le x^{\kappa_0}$.

\medskip
\noindent
  {\bf Triangles with $a_2,~b_2$ and $d$ small.}
\par\noindent We now assume that
$$\max\{a_2,b_2,d\}\le x^{\kappa_0}.
$$
Using \eqref{eq:form1} and
expressing $w$ versus $z$ from the second equation \eqref{eq:form2}
and inserting it into the
first, we get
\begin{equation}
\label{eq:poly}
a_2^2u^4-(-a_2b_2z^2+a_2b_2w^2+a_2d(M^2-N^2))u^2-(b_2wz+dMN)^2=0.
\end{equation}
Fix $d,~a_{2}$ and $b_{2}$. Then $c=d(M^2+N^2)\le x$, therefore we get
$|M|\le (x/d)^{1/2}$ and $|N|\le (x/d)^{1/2}$. In a similar vein,
$b=b_{2}(z^2+w^2)$, therefore $|w|\le (x/b_{2})^{1/2}$, and
$|z|\le (x/b_{2})^{1/2}$. Finally, $|u|\le (x/a_{2})^{1/2}$.
Thus, if the triple $(d,a_{2},b_{2})$ is fixed,
then the parameters $(M,N,w)$ can be fixed in
$O(x^{3/2}/(d{\sqrt{b_{2}}}))$ ways. Now equation \eqref{eq:poly} becomes
an equation of the type $f(u,z)=0$ in the variables $u$ and $z$ alone,
where $f(t,x)\in \Q[t,x]$ is a polynomial of the form
appearing in Lemma \ref{lem2} with
$$
A=a_2^2,\quad B=-a_2b_2,\quad C=a_2b_2w^2+d(M^2-N^2),\quad D=b_2w,\quad
E=dMN.
$$
Furthermore, $u$ and $z$ satisfy  $\max\{|u|,|z|\}\le T$,
where we can take $T=(x/a_{2})^{1/2}$. Since obviously $ABD\ne 0$ for our
polynomial, by Lemma \ref{lem2},
it follows that we can  take $\delta=1$ in Lemma \ref{lem1}, and get
that the number of such solutions for a fixed choice of coefficients is
$$
\ll \left(\frac{x}{a_{2}}\right)^{1/3}.
$$
Hence, for fixed $d,~a_{2}$ and $b_{2}$ the number of choices is
$$
\ll\frac{x^{11/6}}{db^{1/2}_{2}a^{1/3}_{2}}.
$$
Summing the above relations up for $d\le x,~a_2\le x^{\kappa_0}$
and $b_2\le x^{\kappa_0}$,
we get that the number of Heron triangles in this case is
\begin{eqnarray*}
\label{eq:last}
& \ll & x^{11/6}\sum_{b_2<x^{\kappa_0}}\frac{1}{b_2^{1/2}}
\sum_{a_2<x^{\kappa_0}}\frac{1}{a_2^{1/3}}\\
& \ll & x^{11/6}\int_1^{x^{\kappa_0}}\frac{dt}{t^{1/2}}\int_{1}^{x^{\kappa_0}}
\frac{dt}{t^{1/3}}\\
& \ll & x^{11/6}\cdot x^{(1/2+2/3)\kappa_0}=x^{11/6+7\kappa_0/6}.
\end{eqnarray*}
Optimizing, we get $11/6+7\kappa_0/6=2-\kappa_0$, therefore $\kappa_0=1/13$,
which leads to the desired conclusion.
\end{proof}

\noindent
{\bf Acknowledgements.}

\end{document}